\documentclass{amsart}
\usepackage{latexsym}
\usepackage{epsf}
\usepackage{amsfonts}
\usepackage{amsmath}
\usepackage{stmaryrd}
 \newcommand{\resp}{{\it resp.} }
\newcommand{\cf}{{\it cf.} }

\newcommand{\Q}{\mathbb{Q}}
\newcommand{\R}{\mathbb{R}}
\newcommand{\C}{\mathbb{C}}
  \newcommand{\Z}{\mathbb{Z}}

\newcounter{spec}

\swapnumbers

\newtheorem{thm}{Theorem}[subsection]

\theoremstyle{definition}

\newtheorem{ex}[thm]{Example}
\newtheorem{exs}[thm]{Examples}

\begin{document}
\begin{sloppypar} 
 
\title{\bf Ambiguity Theory, old and new}

\author{Yves
Andr\'e}
 
  \address{D\'epartement de Math\'ematiques et Applications, \'Ecole Normale Sup\'erieure  \\ 
45 rue d'Ulm,  75230
  Paris Cedex 05\\France.}
\email{andre@dma.ens.fr}
%  \date{1/04/2007}
  \begin{abstract}  This is a introductory survey of some recent developments of ``Galois ideas" in Arithmetic, Complex Analysis, Transcendental Number Theory and Quantum Field Theory, and of some of their interrelations.
  \end{abstract}
\maketitle

 \bigskip

\section{Introduction. Ambiguity theory according to Galois.}

\subsection{} {\it Ambiguity theory} was the name which Galois used when he referred to his own theory and its future developments. For instance, in his celebrated last letter (1832), he wrote:

\begin{quotation} ``Mes principales m\'editations depuis quelque temps \'etaient dirig\'ees sur l'application \`a l'analyse transcendante de la th\'eorie de l'ambi\-gu\"{\i}t\'e.

 Il s'agissait de voir {\it a priori} dans une relation entre quantit\'es ou fonctions transcendantes quels \'echanges on pouvait faire, quelles quantit\'es on pouvait substituer aux quantit\'es donn\'ees sans que la relation p\^ut cesser d'avoir lieu. Cela fait reconna\^{\i}tre tout de suite l'impossibilit\'e de beaucoup d'expressions que l'on pourrait cher\-cher. 
 
 Mais je n'ai pas le temps et mes id\'ees ne sont pas encore bien d\'evelopp\'ees sur ce terrain qui est immense..." \end{quotation}

\medskip The aim of the present text is to briefly survey the fulfilment of these ideas in several areas of contemporary mathematics. 

More precisely, we shall outline several instances of a common theme, namely: how ambiguities which sometimes occur in the solutions of a well-posed problem, and are first felt as a nuisance, turn out to manifest a rich hidden structure and to provide a decisive tool in  understanding the original problem.

\subsection{} Let us begin with  Galois's own definition of the ``Galois group" attached to an algebraic equation: 

\begin{quotation}
``Soit une \'equation donn\'ee, dont $a,b,c,\dots, $ sont les $m$ racines. Il y aura toujours un {\it groupe de permutations} des lettres $a,b,c,\dots, $ qui jouira de la propri\'et\'e suivante:

$1)$ que toute fonction des racines, invariante par les substitutions de ce groupe, soit rationnellement connue\footnote{Galois adds the following explanation: \begin{quotation}
``Quand nous disons qu'une fonction est {\it rationnellement connue}, nous voulons dire que sa valeur num\'erique est exprimable en fonction rationnelle des coefficients de l'\'equation et des quantit\'es adjointes."\end{quotation}};

$2)$ r\'eciproquement, que toute fonction des ra\-ci\-nes, d\'etermin\'ee rationnellement, soit invariante par ces substitutions."\end{quotation}

 This is probably the first explicit apparition of the notion of {\it group} in mathematics.  

\medskip In modern terms, the set of ``rationally known quantities" forms a field $K$, the set obtained by adjoining to $K$ all polynomial expressions in the roots of the equation is a field extension  $L= K[a,b, c, \dots ]$, and the Galois group, which is defined as 
\begin{equation}\label{Galois} Gal(L/K) := Aut_K \,L,\end{equation}
 is a subgroup of the permutation group of the roots $a,b,c, \ldots $

 Using this concept, and the notions of normal subgroup and solvable group, Galois  proved the impossibility of solving general algebraic equations of degree $\geq 5$ by radicals. After Abel's (slightly earlier) breakthrough, this represents the achievement  of Lagrange's insights about the role of permutations in the reduction to equations of lower order via ``Lagrange resolvents." 

  \begin{quotation}
``Aussi nous savons, nous, ce que cherchait \`a deviner Lagrange, quand il parlait de m\'etaphy\-si\-que \`a propos de ses travaux d'alg\`ebre; c'est la th\'eorie de Galois, qu'il touche presque du doigt, \`a travers un \'ecran qu'il n'arrive pas \`a percer. L\`a o\`u Lagrange voyait des analogies, nous voyons des th\'eor\`emes. Mais ceux-ci ne peuvent s'\'enoncer qu'au moyen de notions et de ``structures" qui pour Lagrange n'\'etaient pas encore des objets math\'e\-ma\-ti\-ques..." \cite{We} \end{quotation}

Eloquent as it is, this famous comment by A. Weil hides the fact that ``Galois theory" is the fruit of a long collective process of multiple reelaborations, lasting until the beginning of the XXth century. We refer to C. Ehrhardt's  recent thesis \cite{E} for a detailed analysis of the mathematical context of Galois's original work and of the slow and indirect construction of ``Galois theory" from Galois to Artin.

 \subsection{}  Although one should be careful, from the viewpoint of the social history of mathematics, when attributing to Galois himself the full conception of the powerful and radical ideas which irrigate ``Galois theory", it is nevertheless worth underlining the thrust of these ideas in modern mathematics:
 
 1)  the emergence of a corpus of abstract concepts or ``structures": groups, normal subgroups, \ldots
 
 2) the emphasis on operations, as abstracted from their outcome; according to the visionary teenager,
 \begin{quotation} ``Sauter \`a pieds joints sur les calculs, grouper les op\'erations, les classer suivant leurs difficult\'es et non leur forme, telle est suivant moi la mission des g\'eom\`etres futurs."  \end{quotation} 
 
 3) the fecundity of the principle of Galois correspondence in various domains (in geometry, it inspired F. Klein's Erlangen program); according to S. Lie,
 \begin{quotation} ``la grande port\'ee de l'oeuvre de Galois tient au fait que sa th\'eorie originale des \'equations alg\'ebriques est une application syst\'ematique des deux notions fondamentales de {\it groupe} et d'{\it invariant}, notions qui tendent \`a dominer la science math\'ematique."  \end{quotation} 
 
 4) the idea of turning obstructions into mathematical objects: do ambiguities form a nuisance? No, they form a group!

\subsection{} Let us end this introduction by mentioning, as a side remark, that some forty years after Galois, the applied mathematician Boussinesq took a completely different attitude in the face of ambiguities (in the context of differential equations of classical mechanics):

\begin{quotation} 
``On trouvait dans le M\'emoire de Poisson des probl\`emes tr\`es simples de M\'ecanique du point mat\'eriel pour lesquels l'\'equation diff\'e\-ren\-tielle permet {\it deux} solutions possibles pour la suite du mouvement du point partant d'une position et d'une vitesse d\'etermi\-n\'ee. C'\'etait, de l'avis de Poisson, une sorte de paradoxe [...]. Boussinesq voit l\`a une {\it propri\'et\'e myst\'erieuse} des solutions singuli\`eres. Ces ind\'etermina\-tions lui paraissent laisser place \`a l'intervention d'un pouvoir directeur distinct des forces m\'ecaniques, c'est-\`a-dire des causes modifiant les acc\'el\'erations. Il voit l\`a une conciliation possible du d\'eterminisme physico-chimique avec l'existence de la vie et de la libert\'e morale." \cite[p. 110]{P} \end{quotation}

\section{Arithmetic versus geometric absolute Galois groups.}

\subsection{The absolute Galois group of $\mathbb Q$.}

A normal extension $K/k$ of fields being fixed, the Galois correspondence relates intermediate field extensions $k \subset \ell \subset K$  and subgroups $Gal(K/ \ell)$ of the Galois group $Gal(K/ k)$. 

 One can also fix only the base field $k$, say  $k={\mathbb Q}$, and let the normal extension $K/{\mathbb Q}$ vary (inside the field $\mathbb C$ of complex numbers, say); one thus obtains a projective system of finite Galois groups $$ \cdots \to Gal(K/{\mathbb Q})  \to \cdots \to  Gal({\mathbb Q}/{\mathbb Q})= \{1\}  .$$
The limit   $$Gal({\bar{\mathbb Q}}/{\mathbb Q}) := \displaystyle\varprojlim_K\, Gal(K/{\mathbb Q})\,$$
is the compact (profinite) group of automorphisms of the field ${\bar{\mathbb Q}}$ of all algebraic numbers.  
  
\smallskip This is the process which the philosopher A. Lautman called ``la mont\'ee vers l'absolu" \cite[III]{Lau}: a single mathematical object, $Gal({\bar{\mathbb Q}}/{\mathbb Q}) $, the absolute Galois group of ${\mathbb Q}$, encodes the properties of all algebraic equations with rational coefficients at the same time.  
 However, the structure of this central object of Number Theory still remains largely mysterious.

\subsection{Galois theory of coverings.} Klein's 1877 remark, that the isometry group of the icosahedron is isomorphic to the Galois group of an algebraic equation of degree $5$ with coefficients in the field of rational functions, opened the door to a ``geometric twin" of Galois theory for coverings (of Riemann surfaces, or more general manifolds). 

In this context, one can again perform ``la mont\'ee vers l'absolu": this leads to the concept of universal covering and of the Poincar\'e fundamental group $\pi_1(X, \ast)$. Following A. Grothendieck, one can get a unified construction (of the arithmetic and geometric versions), at the cost of replacing $\pi_1(X, \ast)$ by its profinite completion\footnote{{\it i.e.} the inverse limit of its finite quotients.} $\hat \pi_1(X, \ast)$: this completion can then be interpreted as the group of automorphisms of the fiber functor 
$$Y \mapsto Y_ \ast $$
on the category of finite coverings of $X$ ($Y_ \ast$ stands for the finite set of points of $Y$ lying above the point $ \ast \in X$). 

\subsection{A geometric approach to $Gal({\bar{\mathbb Q}}/{\mathbb Q}) $.} Let now $X$ be complex plane  deprived of the points $0,1$, so that $\pi_1(X)$ is a free group on the loops $\gamma_0$ and $\gamma_1$ around $0$ and $1$ respectively. 

According to Grothendieck and Belyi, the absolute arithmetic Galois group $Gal({\bar{\mathbb Q}}/{\mathbb Q}) $ acts\footnote{the action depends on the choice of a rational base point $ \ast$, for instance $1/2$, or better, a tangential base point; but the outer action is canonical.}  on the absolute geometric Galois group  $\hat \pi_1(X)$, and this action is faithful: $Gal({\bar{\mathbb Q}}/{\mathbb Q}) $ embeds into the group of (outer) automorphisms $Out\, \hat \pi_1(X)$ of the profinite free group on $\gamma_0$ and $\gamma_1$. In Grothendieck's program, $X = {\mathbb C} \setminus \{0,1\}$ is viewed as the moduli space of four points on the projective line, and $Gal({\bar{\mathbb Q}}/{\mathbb Q}) $ should be studied via its action on the (Grothendieck-Teichm\"uller) tower of configuration spaces of any number of points on the projective line (and on higher genus curves) as well.

\smallskip  A decisive progress occurred around 1990 with V. Drinfeld's study - in a surprising connection with quantum groups - of the {\it Grothendieck-Teichm\"uller group},  a profinite group defined by explicit generators\footnote{consisting in a unit in the profinite free group on one generator together with an element  of the derived group of the profinite free group on two generators.} and relations: two ``hexagonal" and one ``pentagonal" relation (recently, H. Furusho \cite{F} showed that the pentagonal relation\footnote{which comes from the consideration of the configuration of five points on the projective line.} actually implies the hexagonal ones).
 
The Grothendieck-Teichm\"uller group  provides a kind of geometric ``upper bound" for $Gal({\bar{\mathbb Q}}/{\mathbb Q}) $, and it is still unknown whether the two coincide (\cf \cite{LS}).

\section{ Galois ambiguities in ``transcendental analysis". Differential Galois theory.}

Let us now turn to some recent applications of Galois's ideas in ``transcendental analysis", in the spirit of Galois's lines quoted at the beginning of this paper.

\subsection{The beginnings.} Such {\it ambiguities} already puzzled Gauss in his work on the arithmetico-geometric mean and hypergeometric functions. 
The Gauss hypergeometric function
$$F(a,b,c; x)= 1+ \frac{a.b}{1.c}x+ \frac{a(a+1).b(b+1)}{1.2.c(c+1)}x^2+\cdots$$
is a solution of the hypergeometric differential equation 
$$ \;\;\;\;\; x(1-x)F''+ (c-(a+b+1)x)F' - abF=0.$$
The arithmetico-geometric mean is related to $F(\frac{1}{2}, \frac{1}{2}, 1; x)$.
 On the other hand, $F(\frac{1}{4}, \frac{1}{4}, 1; 4x(1-x))$ converges inside the lemniscate
$$\vert 4x(1-x)\vert <1$$ and satisfies there the same hypergeometric differential equation with parameters $a=b=1/2, c=1$.
 Actually $$F(\frac{1}{4}, \frac{1}{4}, 1; 4x(1-x))= F(\frac{1}{2}, \frac{1}{2}, 1;  x)$$ in the left region of the lemniscate, whereas $$F(\frac{1}{4}, \frac{1}{4}, 1; 4x(1-x))= F(\frac{1}{2}, \frac{1}{2}, 1;  1-x)$$ in the right region. 
 
 \smallskip 
  This kind of problems was to be elucidated by Riemann's study \cite{R} of monodromy (that is, in essence, the study of the action of $\pi_1({\mathbb C}\setminus\{0,1\}, \ast)$ on the solutions of the hypergeometric equation).  
  
  \smallskip\noindent  In the same vein, one can also mention the following classical works ({\it cf.} \cite{Gra} for a detailed account of this story):
 
   \smallskip\noindent 
$1)$ Schwarz' classification of hypergeometric differential equations with algebraic solutions:   by looking at the the Galois groups of the corresponding Galois extensions of $\C(x)$, Schwarz reduces it to the classification of finite subgroups of $GL_2$ \cite{Schw}, 
 
   \smallskip\noindent 
$2)$ Jordan's celebrated theorem on finite subgroups of $GL_m$, which arose, via Galois's philosophy, from the study of algebraic solutions of differential equations of order $m$ \cite{J},

  \smallskip\noindent 
$3)$ much of Klein's work, which he described himself as a try to ``blend Galois with Riemann".

 \medskip
\subsection{ Differential Galois theory.} 
This theory, initiated by Picard, Vessiot, Drach and later Kolchin, realizes Galois's program\footnote{a notable precursor of differential Galois theory was Liouville, who not only exhumed and advertised Galois's work but also made the first steps towards Galois's dreams by studying the conditions under which a linear differential equation can be solved by quadratures.}:

\begin{quotation}``voir {\it a priori} dans une relation entre quantit\'es ou fonctions transcendantes quels \'echanges on pouvait faire, quelles quantit\'es on pouvait substituer aux quantit\'es donn\'ees sans que la relation p\^ut cesser d'avoir lieu"\end{quotation}

\noindent when the relation is a {\it linear differential equation $\mathcal L$ of order $m$}. In that case, the set $Sol({\mathcal L}, \ast)$ of solutions of $\mathcal L$ (at some point $\ast$) form a $m$-dimensional vector space over $\C$, the set of ``rationally known quantities" forms a differential field $(K, \partial )$, the ring extension $L $ obtained by adjoining all solutions of $\mathcal L$ and their derivative is a differential field. In practice, $K=\C(x)$ (global case) or $\C((x))$ or $\C(\{x\})\,$ (formal, local case).

\smallskip 
The {\it differential Galois group} 
\begin{equation}\label{diffGalois} Gal({\mathcal L})  = Gal_\partial(L/K) := Aut_{(K,\partial)} \,(L,\partial)\end{equation}
 is an {\it algebraic subgroup} of the group $Aut(Sol({\mathcal L}, \ast)) \cong GL_m(\C)$ of invertible linear transformations of the solution space (it is worth underlining that some basic facts about linear algebraic groups were established by Kolchin for the needs of differential Galois theory).
 
  It turns out that the transcendance degree of $L/K$ (= the maximal number of $K$-algebraically independent solutions of $\mathcal L$) is 
\begin{equation}\label{trdeg}{{\rm tr.deg}_K\, L  = \dim \, Gal({\mathcal L}) .}\end{equation}

\begin{ex}\label{ex321} $K$ as above, $ \;\partial= \frac{d}{dx},\;\;{\mathcal L}:\, x^2\partial y  + y=0$. 

The Galois ambiguities are the ``scaling" $\, e^{\frac{1}{x}}\mapsto \mu .e^{\frac{1}{x}}.$  The differential Galois group is $\,Gal({\mathcal L})= \C^\times\,$ (a one-dimensional algebraic group, in accordance with the fact that  $e^{\frac{1}{x}}$ is transcendental over $K$).\end{ex}

\subsection{Three types of Galois ambiguities for solutions of linear differential equations.}  
$ $

 \subsubsection{1st type of Galois ambiguity: monodromy.} In the global case ({\it i.e.} over the differential field $K= \C(x)$, more precisely, over the base $X= \C\setminus \{a_1,\ldots, a_n\}$) as well as in the local case ({\it i.e.} over $K= \C(\{x\})$, that is, over the germ $X$ of a punctured disk), one has the monodromy homomorphism:
$$\pi_1(X,\ast)\to Aut(Sol({\mathcal L}, \ast)) \cong GL_m(\C)$$ with image the {\it monodromy group} $Mono({\mathcal L})$ of the linear differential equation $\mathcal L$.

The monodromies along any loop are Galois ambiguities, {\it i.e.}
$$ Mono({\mathcal L}) \subset Gal({\mathcal L})$$

\smallskip Schlesinger's theorem \cite{S} tells that this is the only kind of Galois ambiguities in the Fuchsian case: {\it for any differential equation with regular singularities ({\it i.e.} Fuchsian), $ Mono({\mathcal L}) $ is Zariski-dense in $Gal({\mathcal L})$}.

 \subsubsection{2nd type of Galois ambiguity: exponential scaling:} as in the above example\footnote{which shows that Schlesinger's theorem does not extends in the presence of irregular singularities.} $  e^{\frac{1}{x}}\mapsto \mu .e^{\frac{1}{x}}$.

Over the differential field $K=\C((x))$ ({\it i.e.} in the formal local case), it turns out that formal monodromies and exponential scalings generate a Zariski-dense subgroup of $Gal({\mathcal L})$.

 \medskip\noindent {\it Note}: any linear differential equation over $\C((x))$ has a basis of solutions of the form 
$$\sum (x^{\lambda_i} \log^{k_i} x ). e^{P_i( {x^{-\frac{1}{e}}})}. \hat f_i(x)$$ for suitable polynomials $P_i$'s and formal power series $f_i$'s. The exponential scalings act on the factors $e^{P_i( {x^{-\frac{1}{e}}})}$, formal monodromies on the factors $x^{\lambda_i} \log^{k_i} x$.

 \medskip 
 Over $K=\C(\{x\})$ (analytic local case), one encounters a

 \subsubsection{3rd type of Galois ambiguity: Stokes phenomenon.}\footnote{whose discovery, in a special case, dates back to 1857 (same year as the publication of \cite{R}).}
  In general, divergent series $\hat f_i$ occur in the solutions. They may be resummed, using a Borel-Laplace type process, in suitable sectors bissected by an arbitrary fixed half-line $\ell$  (except for finitely many $\ell$'s, the singular directions).

  \smallskip\noindent {\it Example}: $\mathcal L$ $ : \; x^2\partial y +y=x.$
   The formal solution $ \hat y = \sum_0^\infty\, (-1)^nn!x^{n+1}$
  can resummed (in a sector bissected by $\ell= \R_+$) to $\,\int_0^\infty\, \frac{e^{-t/x}}{1+t}dt.$  
   
\smallskip  For a singular line ($\ell = \R_-$ in this example), there are two summations (from above and from below), which are related by a unipotent matrix (Stokes matrix).  These matrices are ``Galois transformations", {\it i.e.} belong to $Gal({\mathcal L})$.

\smallskip Ramis' theorem tells that, in the local case ($K=\C(\{x\})$) as well as in the global case ($K=\C(x)$),  those three kinds of Galois ambiguities are the only ones: {\it formal monodromies, exponential scalings and Stokes matrices generate a Zariski-dense subgroup of $Gal({\mathcal L})$}.

More precisely, on  formalizing the algebraic relations between these explicit Galois transformations, Ramis has determined the structure of the tannakian group of the category of linear differential equation over $\C(\{x\})$, 
 the {\it wild fundamental group} $\pi_1^{wild}$, \cf \cite{vdPS} for a detailed account of this story.

\smallskip
{\it Example of application}: J. Morales and J.P. Ramis \cite{MR} have given a Galois criterion for the {\it Liouville integrability of certain hamiltonian systems} in celestial mechanics (essentially, the differential Galois group of the linearized equation should be virtually abelian); thanks to the above description of topological generators of the Galois group, this criterion could be checked in a number of concrete instances.  

\subsection{ Galois ambiguities for non-linear dynamical systems? }
 Actually, the early prospects of differential Galois theory, especially in Drach's plans, were much more ambitious. But in spite of the successes of differential algebra (Ritt, ...), it soon appeared that the foundations of differential Galois theory were not solid enough outside of the linear case.  
 
 After various attempts by J.F. Pommaret (using Lie pseudo-groups) and H. Umemura (using infinite-dimensional formal groups), a new theory built by B. Malgrange seems extremely promising \cite{M}. In the beginning, Malgrange\footnote{and P. Cartier, independently.} remarked that another way to ``save" Schlesinger's theorem $$ \overline{Mono({\mathcal L})}^{Zar} = Gal({\mathcal L})$$ for non-Fuchsian linear differential equations - without adding new Galois ambiguities - is to {\it replace algebraic groups by algebraic groupoids}.
 This approach then generalizes to the non-linear case (foliations with singularities) if one further replace "algebraic groupoid" (defined by algebraic equations) by ``algebraic $D$-groupoid" (defined by algebraic systems of partial differential equations, using jets).

 The Galois $D$-groupoid of an algebraic or analytic foliation is (essentially) defined as the {\it smallest $D$-groupoid of the groupoid of germs of invertible transformations containing the flows of the tangent vector fields}.

\smallskip 
{\it Example of application}: Umemura \cite{U} and G. Casale \cite{Cas} settled the question of irreducibility of Painlev\'e's equation I $$P_I: y''= 6y^2+x,$$ one century after unsatisfactory attempts by Painlev\'e and Drach.

 \bigskip
\section{Galois ambiguities in transcendantal number theory. Motivic Galois theory.}

\subsection{Periods.} Let us now leave functions (and differential equations) and come back to numbers. In this section, we address the question:
does Galois theory extend from its original setting of algebraic numbers to the setting of transcendental numbers?

\smallskip It turns out that the answer is yes, at least conjecturally, if one restricts oneself to so-called {\it periods}; we shall briefly explain how to attach to a period its conjugates, and a Galois group which permutes them transitively.

\smallskip A {\it period} is a complex number whose real and imaginary parts are absolutely convergent multiple integrals 
   $$\alpha = \frac{1}{{(2\pi )}^m}\, \int_\Sigma \,\Omega$$ where $\Sigma $ is a domain in $\R^n$ defined by polynomial equations and inequations with rational coefficients, $\Omega$ is a rational differential form with rational coefficients, and $m$ is a natural integer. The set of periods is a countable subring of $\C$ which contains $\bar\Q$.
   
 \begin{exs} $1)$ The very name ``period" comes from the case of elliptic periods (in the case of an elliptic curve defined over $\bar\Q$, the periods in the classical sense are indeed periods in the above sense). 
 
 \noindent $2)$ The values at algebraic numbers $x$ of the hypergeometric function $F(a,b,c;x)$ (and of any generalized hypergeometric functions ${}_pF_{p-1}$) with rational parameters are periods.
  
  \noindent $3)$  Periods also frequently appear in quantum field theory: P. Belkale and P. Brosnan \cite{BB} show that Feynman amplitudes $I({\rm D})$ with rational parameters 
 can be written as a product of a Gamma-factor and a meromorphic function $H({\rm D})$ such that the coefficients of its Taylor expansion at any integral value of ${\rm D}$ are all periods. \end{exs}
 
  Since periods are on one hand transcendental numbers in general, and on the other hand are related to algebro-geometric objects, one may expect that the corresponding expected Galois groups will not be finite groups, but linear algebraic groups (just as in the case of transcendental solutions of linear differential equations with rational coefficients). We will see that it is essentially so (at least conjecturally) and that the answer belongs to the framework of motivic Galois theory\footnote{motivic Galois theory was first envisioned by Grothendieck in the mid sixties.} (for a more detailed account of the story, we refer to \cite{A2}).
  
 \subsection{Motives.}
Let $X$ be a smooth algebraic variety over $\Q$, and let $Y$ be a closed (possibly reducible) subvariety. Periods arise as entries of a matrix of the comparison isomorphism, 
 given by integration of algebraic diffential forms over chains, between 
 algebraic De Rham and ordinary Betti relative cohomology:
   \begin{equation}\label{isocomp} H_{DR}(X,Y)\otimes \C \stackrel{\varpi_{X,Y}}{\cong}  H_{B}(X,Y)\otimes \C.\end{equation}
   This is where motives enter the stage.   Motives are intermediate between algebraic varieties and their linear invariants (cohomology): they are of algebro-geometric nature on one hand, but they are supposed to play the role of a universal cohomology for algebraic varieties and thus to enjoy the same formalism on the other hand. 
    One expects the existence of an abelian category ${\rm MM} $ of {\it mixed motives} (over $\Q$, with rational coefficients), and of a functor
    $$h:  Var(\Q) \to {\rm MM}$$ (from the category of algebraic varieties over $\Q$) which plays the role of universal cohomology. The morphisms in ${\rm MM}$ should be related to algebraic correspondences.
     In addition, the cartesian product on $Var(\Q)$ corresponds via $h$ to a certain tensor product $\otimes$ on ${\rm MM}$, which makes ${\rm MM}$ into a {\it tannakian category}, {\it i.e.} it has the same formal properties as the category of representations of a group. 
 
   \smallskip   The cohomologies $H_{DR}$ and $H_B$ factor through $h$, giving rise to two $\otimes$-functors 
 $$ H_{DR}, \, H_B:\, {\rm MM} \to Vec_\Q$$ with values in the category of finite-dimensional $\Q$-vector spaces. 
   Moreover, corresponding to \eqref{isocomp},  there is a isomorphism in $Vec_\C$
 \begin{equation}  \varpi_M:  H_{DR}(M)\otimes \C \cong H_B(M)\otimes \C  \end{equation} which is
 $\otimes$-functorial in the motive $M$. The entries of a matrix of $\varpi_{M}$ with respect to some basis of the $\Q$-vector space $H_{DR} (M)$ (\resp $H_{B} (M )$) are the {\it periods} of $M$.
          
      \subsection{Motivic Galois groups and period torsors.}  
  Let $\langle M \rangle$ be the tannakian subcategory of ${\rm MM}$ generated by a motive $M$: its objets are given by algebraic constructions on $M$ (sums, subquotients, duals, tensor products).
  
  One defines the {\it motivic Galois group} of $M$ to be the group scheme
  \begin{equation} G_{mot}(M) := Aut^\otimes \,{H_B}_{\mid \langle M \rangle}  \end{equation}
  of automorphisms of the restriction of the $\otimes$-functor ${H_B}$ to $  \langle M \rangle$. 
  
This is a linear algebraic group over $\Q$: in heuristic terms, $ G_{mot}(M)$  is just {\it the Zariski-closed subgroup of $GL ( H_B(M))$  consisting of matrices which preserve motivic relations in the algebraic constructions on $H_B(M)$}.
  
  \smallskip Similarly, one can consider both $H_{DR}$ and $H_B$, and define the {\it period torsor} of $M$ to be 
      \begin{equation} P_{mot}(M) := Isom^\otimes \,({H_{DR}}_{\mid \langle M \rangle}, {H_B}_{\mid \langle M \rangle}) \in Var(\Q)  \end{equation}
  of isomorphisms of the restrictions of the $\otimes$-functors ${H_{DR}}$ and ${H_B}$ to $  \langle M \rangle$. 
This is a torsor under  $ G_{mot}(M) $, and it has a canonical complex point:
   \begin{equation}\varpi_M\in  P_{mot}(M)(\C).   \end{equation}

 \subsection{Grothendieck's period conjecture.} It asserts that $\varpi_M$ is a \emph{generic point}: 

{\it the smallest algebraic subvariety of $P_{mot}(M)$ defined over $\Q$ and containing $\varpi$ is $P_{mot}(M)$ itself. }
  
  \medskip 
 In more heuristic terms, this means that any polynomial relations with rational coefficients between periods should be of motivic origin (the relations of motivic origin being precisely those which define     $P_{mot}(M)$).

    \smallskip The conjecture is also equivalent to: 
 {\it  $P_{mot}(M)$ is connected (over $\Q$) and   
  \begin{equation}\label{trdegmot}{\rm tr. \,deg}_\Q \,  \Q[{\rm periods}(M)] = \dim \,G_{mot}(M) . \end{equation} }

 \smallskip {\it Example:} for the motive of ${\mathbb P}^1$, the motivic Galois group is the mutiplicative group ${\mathbb G}_m$, a period is $ 2\pi i$, and Grothendieck's conjecture amounts to the transcendence of $\pi$. 
 
  \smallskip {\it Remarks.} $1)$ There is an analogy between equations \eqref{trdeg} and \eqref{trdegmot}, which can be more made more precise if one considers a family of smooth algebraic varieties (or motives): the variation of periods is then controlled by a linear differential equation (the so-called Picard-Fuchs equation. For the relation between its differential Galois group and the motivic Galois groups of the fibers, we refer to \cite[\S 6]{A2}. 
  
 \smallskip \noindent $2)$ By definition, periods are convergent integrals of a certain type. They can be transformed by algebraic changes of variable, or using additivity of the integral, or using Stokes formula. 
  
    Kontsevich has conjectured that {\it any polynomial relation with rational coefficients between periods can be obtained by way of these elementary operations from calculus} (\cf \cite{KZ}). Using ideas of Nori and the expected equivalence of various motivic settings, it can be shown that this conjecture is actually equivalent to Grothendieck's conjecture (\cf \cite[ch. 23]{A1}).

   \subsection{Outline of a Galois theory for periods.}  Let $\alpha$ be a period. There exists a motive $M\in {\rm MM}$ such that $\alpha\in  \Q[{\rm periods}(M)]$.
   Let us assume Grothendieck's period conjecture for $M$. Then $ \Q[{\rm periods}(M)] $ coincides with the algebra $\Q[P_{mot}(M)]$ of functions on $P_{mot}(M)$. Since $P_{mot}(M)$ is a torsor under $G_{mot}(M)$, the group of rational points $G_{mot}(M)(\Q)$ acts on $\Q[P_{mot}(M)]$, hence on $ \Q[{\rm periods}(M)] $. 
   
\smallskip   One defines the {\it conjugates of $\alpha$} to be the elements of the orbit $G_{mot}(M)(\Q).\alpha $. It follows from Grothendieck's conjecture that this does not depend on $M$. 
   
 \smallskip   The {\it Galois closure} $\Q[\alpha]_{gal}$ of $\Q[\alpha]$ is the subalgebra $\Q[G_{mot}(M)(\Q).\alpha ]$ of $\Q[{\rm periods}(M)]$.
   
 \smallskip   The {\it Galois group} of $\alpha$ is the smallest quotient  $G_\alpha$ of $ G_{mot}(M)(\Q)$ which acts effectively on $\Q[\alpha]_{gal}$. It acts transitively on the set of its conjugates.

\subsection{Examples.} $1)$ The conjugates of $2\pi i$ (a period of $M=h({\mathbb P}^1)$) are its non-zero rational multiples, $G_\pi = \Q^\times$.

\medskip\noindent $2)$ Let $x$ be a non-zero rational number of absolute value than $1$. Then $\alpha = F(\frac{1}{2}, \frac{1}{2}, 1;  x)$ is, up to multiplication by $2\pi i$, a period of the elliptic curve $ Y^2= X(X-1)(X-x)$. 
 Except for finitely many exceptions\footnote{the so-called singular moduli, corresponding to elliptic curves with complex multiplication.} $x\in \Q$, 
 $G_\alpha = GL_2(\Q)$, and 
 the conjugates  of $\alpha$ are
    $ mF(\frac{1}{2}, \frac{1}{2}, 1;  x) + 
     n {F( \frac{1}{2}, \frac{1}{2}, 1;  {1-x})},\; m,n\in \Q$, not both $0$. 
 
\medskip\noindent $3)$ Let $s$ be an odd integer $>1$. Then $ \zeta(s):= \sum n^{-s}$ is a period of a so-called mixed Tate motive over $\Z$. Grothendieck's period conjecture for this type of motives would imply that $\pi$ and $\zeta(3),\zeta(5),\ldots$ are algebrically independent, that the conjugates of $\zeta(s)$ are  $\zeta(s)+ r (\pi i)^s, \;r\in \Q,$ and that $G_{\zeta(s)}$ is a semi-direct product of $\Q^\times$ by $\Q$.

\medskip\noindent $4)$ More generally, multiple zeta values  $$\displaystyle{  \zeta (s_1,\ldots, s_k) = \sum_{n_1>\cdots> n_k>0} \, \frac{1}{n_1^s\ldots n_k^{s_k}} }  
 = \int_{1>t_1>\cdots >t_s>0} \, \omega_{s_1}\wedge \cdots \wedge \omega_{s_k},$$ 
  (where  $\;\; s_i>0, s_1>1, s= s_1+\cdots +s_k, \, \omega_0= \frac{dt}{t}, \omega_1=\frac{dt}{1-t}, \omega_r= \omega_0^{r-1}\wedge \omega_1$ for $r>1$) are also periods for mixed Tate motives over $\Z$.
 Many polynomial relations are known among these numbers (double shuffle and regularization relations,  Drinfeld's associator relations);
so far, all these relations have been proven to be of motivic origin, as predicted by Grothendieck's period conjecture. 

The ``absolute" motivic Galois group corresponding to multiple zeta values (taken altogether) is expected to be  
   \begin{equation} {U^\ast(3,5,7\ldots)= U(3,5,7\ldots)\propto {\mathbb G}_m}  \end{equation}  
   where 
   $U(3,5,7\ldots)$ is prounipotent group whose Lie algebra, graded by the ${\mathbb G}_m$-action, is the free graded Lie algebra with one generator in each odd degree $3,5,7\ldots$. 

\smallskip The fact that multiple zeta values satisfy Drinfeld's associator relations provides a homomorphism 
 \begin{equation} U^\ast(3,5,7\ldots) \to GT \end{equation} where $GT$ denotes the ``proalgebraic version" of the (profinite) Grothendieck-Teichm\"uller group mentioned in \S 2.3. It is still unknown whether this is an isomorphism (\cf \cite[ch. 25]{A1}).

\section{Galois ambiguities in quantum field theory.}

In the last years, Galois ideas have barged into quantum field theory. We end this survey by just giving a  glimpse of the main steps of these developments (see \cite{C}, \cite{CM1}). 

\subsection{} D. Kreimer, A. Connes \cite{CK}: conceptual understanding of (reg-dim) perturbative renormalization of a QFT in dimension $D$ in terms of Birkhoff decomposition of loops $\gamma(z)$ in a certain pro-unipotent group $U$ attached to the QFT via the construction of a Hopf algebra on suitable Feynman diagrams; here $z$ is the complexified dimension minus $D$.

The postulated independence of $\gamma_-(z)$ from the choice of the auxiliary mass scaling $\mu$ allows to lift the renormalization group to a one-parameter subgroup of $U$ (with infinitesimal generator $\mu\frac{\partial}{\partial \mu}$).

\subsection{} Connes, M. Marcolli \cite{CM2}:  the universal situation is given by $U(1,2,3\ldots)$ and a certain canonical one-parameter subgroup. The Birkhoff decomposition can be reinterpreted in terms of a Riemann-Hilbert correspondence for a category of so-called ``equisingular" integrable connections over 
$$\{0<\vert z\vert <1\}\times \{\mu \neq 0\}.$$ The corresponding universal differential Galois group is  
 \begin{equation} {U^\ast(1,2,3\ldots)= U(1,2,3\ldots)\propto {\mathbb G}_m} \end{equation}    (Cartier's {\it cosmic Galois group} \cite{Ca}).

Hence perturbative renormalization is given a differential Galois interpretation, where the Galois ambiguities are of type 1 (monodromies corresponding to the loops $\gamma_-(z)$) or 2 (``exponential scaling" corresponding to the mass scaling $\mu$). There is no Stokes phenomenon here.

  \subsection{} On the other hand, $ {U^\ast(1,2,3\ldots) }$ is abstractly isomorphic to the universal motivic Galois group attached to mixed Tate motives over ${\mathbb Z}[i, \frac{1}{2}]$. 
    
    To complete the picture, it remains to understand precisely the role of  mixed Tate motives in perturbative renormalization; in this direction, S. Bloch, H. Esnault and Kreimer \cite{BEK} have given a motivic interpretation of the Kreimer-Connes construction. 

\bigskip In conclusion, one may say that the ambiguities arising from divergences in QFT, far from being a nuisance, are the sign of the existence of a rich and mysterious structure with many Galois symetries.

 \end{sloppypar}
 \end{document}